\documentclass[10pt,twoside]{amsart}

\newtheorem{theorem}{theorem}[section]
\newtheorem{lemma}[theorem]{Lemma}
\newtheorem{prop}[theorem]{Proposition}
\newtheorem{cor}[theorem]{Corollary}

\newtheorem{remark}[theorem]{Remark}

\def\YEAR{\year}\newcount\VOL\VOL=\YEAR\advance\VOL by-1995
\def\firstpage{1}\def\lastpage{1000}
\def\received{}\def\revised{}
\def\communicated{}

\def\a{{\mathfrak a}}

\def\c{{\mathfrak c}}
\def\d{{\mathfrak d}}
\def\g{{\mathfrak g}}
\def\h{{\mathfrak h}}
\def\k{{\mathfrak k}}
\def\l{{\mathfrak l}}
\def\n{{\mathfrak n}}

\def\p{{\mathfrak p}}
\def\q{{\mathfrak q}}

\def\t{{\mathfrak t}}
\def\uu{{\mathfrak u}}

\def\z{{\mathfrak z}}

\def\R{{\mathbb R}}

\makeatletter
\def\magnification{\afterassignment\m@g\count@}
\def\m@g{\mag=\count@\hsize6.5truein\vsize8.9truein\dimen\footins8truein}
\makeatother

\font\eightrm=cmr8
\font\caps=cmcsc10                    % Theorem, Lemma etc
\font\Caps=cmcsc10 scaled \magstep1   % Title

\makeatletter
\@specialpagefalse\headheight=8.5pt
\def\DocMath{}
\renewcommand{\@evenhead}{%
    \ifnum\thepage>\lastpage\rlap{\thepage}\hfill%
    \else\rlap{\thepage}\slshape\leftmark\hfill{\caps\SAuthor}\hfill\fi}%
\renewcommand{\@oddhead}{%
    \ifnum\thepage=\firstpage{\DocMath\hfill\llap{\thepage}}%
    \else{\slshape\rightmark}\hfill{\caps\STitle}\hfill\llap{\thepage}\fi}%
\makeatother

\def\TSkip{\bigskip}
\newbox\TheTitle{\obeylines\gdef\GetTitle #1
\ShortTitle  #2
\SubTitle    #3
\Author      #4
\ShortAuthor #5
\EndTitle
{\setbox\TheTitle=\vbox{\baselineskip=20pt\let\par=\cr\obeylines%
\halign{\centerline{\Caps##}\cr\noalign{\medskip}\cr#1\cr}}%
	\copy\TheTitle\TSkip\TSkip%
\def\next{#2}\ifx\next\empty\gdef\STitle{#1}\else\gdef\STitle{#2}\fi%
\def\next{#3}\ifx\next\empty%
    \else\setbox\TheTitle=\vbox{\baselineskip=20pt\let\par=\cr\obeylines%
    \halign{\centerline{\caps##} #3\cr}}\copy\TheTitle\TSkip\TSkip\fi%
\centerline{\caps #4}\TSkip\TSkip%
\def\next{#5}\ifx\next\empty\gdef\SAuthor{#4}\else\gdef\SAuthor{#5}\fi%
\ifx\received\empty\relax
    \else\centerline{\eightrm Received: \received}\fi%
\ifx\revised\empty\TSkip%
    \else\centerline{\eightrm Revised: \revised}\TSkip\fi%
\ifx\communicated\empty\relax
    \else\centerline{\eightrm Communicated by \communicated}\fi\TSkip\TSkip%
\catcode'015=5}}\def\Title{\obeylines\GetTitle}
\def\Abstract{\begingroup\narrower
    \parskip=\medskipamount\parindent=0pt{\caps Abstract. }}
\def\EndAbstract{\par\endgroup\TSkip}

\long\def\MSC#1\EndMSC{\def\arg{#1}\ifx\arg\empty\relax\else
     {\par\narrower\noindent%
     2000 Mathematics Subject Classification: #1\par}\fi}

\long\def\KEY#1\EndKEY{\def\arg{#1}\ifx\arg\empty\relax\else
	{\par\narrower\noindent Keywords and Phrases: #1\par}\fi\TSkip}

\newbox\TheAdd\def\Addresses{\vfill\copy\TheAdd\vfill
    \ifodd\number\lastpage\vfill\eject\phantom{.}\vfill\eject\fi}
{\obeylines\gdef\GetAddress #1
\Address #2 
\Address #3
\Address #4
\EndAddress
{\def\xs{4.3truecm}\parindent=0pt
\setbox0=\vtop{{\obeylines\hsize=\xs#1\par}}\def\next{#2}
\ifx\next\empty % 1 address
     \setbox\TheAdd=\hbox to\hsize{\hfill\copy0\hfill}
\else\setbox1=\vtop{{\obeylines\hsize=\xs#2\par}}\def\next{#3}
\ifx\next\empty % 2 addresses
     \setbox\TheAdd=\hbox to\hsize{\hfill\copy0\hfill\copy1\hfill}
\else\setbox2=\vtop{{\obeylines\hsize=\xs#3\par}}\def\next{#4}
\ifx\next\empty\ % 3 addresses
     \setbox\TheAdd=\vtop{\hbox to\hsize{\hfill\copy0\hfill\copy1\hfill}
                \vskip20pt\hbox to\hsize{\hfill\copy2\hfill}}
\else\setbox3=\vtop{{\obeylines\hsize=\xs#4\par}}
     \setbox\TheAdd=\vtop{\hbox to\hsize{\hfill\copy0\hfill\copy1\hfill}
	        \vskip20pt\hbox to\hsize{\hfill\copy2\hfill\copy3\hfill}}
\fi\fi\fi\catcode'015=5}}\gdef\Address{\obeylines\GetAddress}

\begin{document}

%\author[P. Foth]{Philip Foth}
%\author[M. Otto]{Michael Otto}
%\address{Department of Mathematics, 
%University of Arizona, Tucson, AZ 85721.}
%\email{foth@math.arizona.edu, otto@math.arizona.edu}
%\title[Symplectic realization of van den Ban's theorem]
%{A symplectic realization of van den Ban's convexity theorem}
%%\thanks{The first author was supported in part by NSF grant DMS-0097314}

%%%%% ------------- fill in your data below this line  -------------------
%%%%%    The following lines \Title ... \EndAddress must ALL be present
%%%%%    and in the given order.
\Title
A symplectic approach to van den Ban's convexity theorem  
%%%%%    Put here the title. Line breaks will be recognized. 
\ShortTitle 
Symplectic approach to van den Ban's theorem
%%%%%    Running title for odd numbered pages, ONE line, please. 
%%%%%    If none is given, \Title will be used instead.          
\SubTitle   
%%%%%    A possible subtitle goes here.
\Author 
Philip Foth, Michael Otto
%%%%%    Put here name(s) of authors. Line breaks will be recognized.  
\ShortAuthor
P. Foth, M. Otto 
%%%%%%   Running title for even numbered pages, ONE line, please. 
%%%%%%   If none is given, \Author will be used instead.          
\EndTitle
\Abstract 
Let $G$ be a complex semisimple Lie group and $\tau$ a complex antilinear 
involution that commutes with a Cartan involution. If $H$ denotes the 
connected subgroup of $\tau$-fixed points in $G$, and $K$ is maximally 
compact, each $H$-orbit in $G/K$ can be equipped with a Poisson structure 
as described by Evens and Lu. We consider symplectic leaves of certain 
such $H$-orbits with a natural Hamiltonian torus action. A symplectic convexity 
theorem then leads to van den Ban's convexity result for (complex) semisimple symmetric spaces. 
%%%%%    Put here the abstract of your manuscript.
%%%%%    Avoid macros and complicated TeX expressions, as this is
%%%%%    automatically translated and posted as an html file.
\EndAbstract
\MSC 
53D17, 53D20, 22E46
%%%%%    2000 Mathematics Subject Classification: 
\EndMSC
\KEY 
Lie group, real form, Poisson manifold, symplectic leaf, moment map, convex cone 
%%%%%    Keywords and Phrases:     
\EndKEY
%%%%%    All 4 \Address lines below must be present. To center the last
%%%%%    entry, no empty lines must be between the following \Address
%%%%%    and \EndAddress lines.
\Address 
Philip Foth \\ Department of Mathematics \\  University of Arizona \\  Tucson, AZ 85721-0089 \\ U.S.A. \\ foth@math.arizona.edu
\Address
Michael Otto \\ Department of Mathematics \\  University of Arizona \\  Tucson, AZ 85721-0089 \\ U.S.A. \\ otto@math.arizona.edu
%%%%%    Address of second Author here etc.
\Address
\Address
\EndAddress
%%

%%       Make sure the last tex command in your manuscript
%%       before the first \end{document} is the command  \Addresses
%%
%%---------------------Here the prologue ends---------------------------------
%%--------------------Here the manuscript starts------------------------------

%\maketitle

\section{Introduction}

In 1982, Atiyah \cite{A} discovered a surprising connection between
results in Lie theory and symplectic geometry. He proved a general 
symplectic convexity theorem of which Kostant's linear convexity 
theorem (for complex semisimple Lie groups) is a corollary. In this
context, the orbits relevant for Kostant's theorem carry the natural
symplectic structure of coadjoint orbits. The symplectic convexity theorem, 
which was found independently by Guillemin and Sternberg \cite{GS}, 
states that the image under the moment map of a compact connected 
symplectic manifold with Hamiltonian torus action is a convex polytope. 
Subsequently, Duistermaat \cite{D} extended the symplectic convexity 
theorem in a way that it could be used to prove Kostant's linear theorem 
for real semisimple Lie groups as well. 

Lu and Ratiu \cite{LR} found a way to put Kostant's nonlinear theorem 
into a symplectic framework. For a complex semisimple Lie group $G$ 
with Iwasawa decomposition $G=NAK$, they regard the relevant $K$-orbit 
as symplectic leaves of the Poisson Lie group $AN$, carrying the 
Lu-Weinstein Poisson structure. Kostant's nonlinear theorem for both 
complex and certain real groups then follows from the AGS-theorem or 
Duistermaat's theorem. 

In this paper, we want to give a symplectic interpretation of van den Ban's 
convexity theorem for a complex semisimple symmetric space $(\g,\tau)$, which is 
a generalization of Kostant's nonlinear theorem for complex groups. For the precise statement 
of van den Ban's result we refer to Section \ref{vdb}. The main difference 
in view of our symplectic approach is that van den Ban's theorem is concerned 
with orbits of a certain subgroup $H\subset G$ that are in general neither 
symplectic nor compact. Since $G$ is complex we can use 
a method due to Evens and Lu \cite{EL} to equip $H$-orbits in $G/K$ with a 
certain Poisson structure. An $H$-orbit foliates into symplectic leaves, 
and on each leaf some torus acts in a 
Hamiltonian way. The corresponding moment map $\Phi$ turns out to be proper, 
and therefore the symplectic convexity theorem of Hilgert-Neeb-Plank \cite{HNP} 
can be applied, which describes the image under $\Phi$ in terms of local moment cones. An analysis of those local moment cones shows that the image of $\Phi$ is the sum of a compact convex polytope and a convex polyhedral cone, just as in van den Ban's theorem.

The case of van den Ban's theorem for a real semisimple symmetric space is dealt with in a separate paper \cite{O}. It follows the symplectic approach of Lu and Ratiu towards Kostant's nonlinear convexity theorem. The main tool is a generalized version of Duistermaat's theorem for non-compact manifolds.
\par
{\it Acknowledgments.} We are grateful to the referee for the careful reading of 
the manuscript and many useful comments and suggestions.

\section{Van den Ban's theorem}\label{vdb}

The purpose of this section is to fix notation and to recall the statement of van den Ban's theorem. 

Let $G$ be a real connected semisimple Lie group with finite center, equipped with an involution $\tau$, 
i.e. $\tau$ is a smooth group homomorphism such that $\tau^2=id$. Let $\g$ be the Lie algebra of $G$. 
We write $H$ for an open subgroup of 
$G^\tau$, the $\tau$-fixed points in $G$. Let $K$ be a $\tau$-stable maximal compact subgroup of $G$. 
The corresponding Cartan involution $\theta$ on $\g$ commutes with $\tau$ and induces the Cartan 
decomposition $\g=\k+\p$. If $\h$ and $\q$ denote the $(+1)$- and $(-1)$-eigenspace of $\g$ 
with respect to $\tau$ one obtains
\[ \g = (\k\cap \h) + (\p\cap \h) + (\k\cap \q) + (\p\cap \q) . \]

We fix a maximal abelian subalgebra $\a^{-\tau}$ of $\p\cap \q$. (In \cite{V} this subalgebra is denoted 
by $\a_{pq}$.) In addition, we choose $\a^{\tau}\subseteq \p\cap \h$ such that $\a := \a^{\tau}+\a^{-\tau}$ 
is maximal abelian in $\p$. Let $\Delta(\g,\a^{-\tau})$ and $\Delta(\g,\a)$ denote the sets of roots for 
the root space decomposition of $\g$ with respect to $\a^{-\tau}$ and $\a$, respectively. Next, we choose  
a system of positive roots $\Delta^+(\g,\a)$ and define 
$$
\Delta^+(\g,\a^{-\tau}) =\{ \alpha |_{\a^{-\tau}} : \alpha\in\Delta^+(\g, \a), \alpha |_{\a^{-\tau}}\ne 0\}.
$$
This leads to an Iwasawa decomposition
\[ \g=\n+\a+\k=\n^1+\n^2+\a+\k , \]
where
\begin{eqnarray*}
\n & = & \sum_{\alpha \in \Delta^+(\g,\a)} \g^{\alpha}, \\
\n^1 & = & \sum_{\alpha \in \Delta^+(\g,\a), \alpha |_{\a^{-\tau}}\neq 0} \g^{\alpha} = \sum_{\beta \in \Delta^+(\g,\a^{-\tau})} \g^{\beta}, \\
\n^2 & = & \sum_{\alpha \in \Delta^+(\g,\a), \alpha |_{\a^{-\tau}}= 0} \g^{\alpha}. 
\end{eqnarray*}
Here $\ \g^{\alpha}=\{ X\in \g : [H,X]=\alpha(H)X \ \forall H\in \a \} $ for $\alpha \in \Delta(\g,\a)$, 
and similarly $\g^\beta$ is defined for $\beta \in \Delta(\g,\a^{-\tau})$.  

Let $N$ and $A$ denote the analytic subgroups of $G$ with Lie algebras $\n$ and $\a$, respectively. 
The Iwasawa decomposition $G=NAK$ on the group level has the middle projection $\mu:G\rightarrow A$. 
We write $pr_{\a^{-\tau}}:\a \rightarrow \a^{-\tau}$ for the projection along $\a^{\tau}$. 

For $\beta \in \Delta^+(\g,\a^{-\tau})$ define $H_\beta \in \a^{-\tau}$ such that
\[ H_\beta \perp \mbox{ker} \beta, \quad \beta(H_\beta)=1 , \]
where $\perp$ means orthogonality with respect to the Killing form $\kappa$. 

Note that the involution $\theta \circ \tau$ leaves each root space 
\[ \g^{\beta}=\sum_{\alpha \in \Delta(\g,\a), \alpha|_{\a^{-\tau}}=\beta} \g^{\alpha} \]
stable. Each $\g^{\beta}=(\g^{\beta})_+ \oplus (\g^{\beta})_- $ decomposes into $(+1)$- and $(-1)$-eigenspace with respect to $\theta \circ \tau$. 

For 
\[ \Delta_{-}:=\{ \beta \in \Delta(\g,\a^{-\tau}) : (\g^{\beta})_- \neq 0 \} , \]
let $\Delta^{+}_{-}=\Delta_{-} \cap \Delta^+(\g,\a^{-\tau})$. Define the closed cone
\[ \Gamma(\Delta^{+}_{-})=\sum_{\beta \in \Delta^{+}_{-}} \R_+ H_\beta . \]
Write $\mathcal{W}_{K \cap H}$ for the Weyl group
\[ \mathcal{W}_{K \cap H}=N_{K \cap H}(\a^{-\tau})/Z_{K \cap H}(\a^{-\tau}) . \]
The convex hull of a Weyl group orbit through $X\in \a^{-\tau}$ will be denoted by 
${\rm conv}(\mathcal{W}_{K \cap H}.X)$. 

\begin{remark}\label{thetatau}
Consider the Lie algebra $\g^{\theta\tau}$ of $\theta\tau$-fixed points in $\g$. It is reductive and its semisimple part $\g'=[\g^{\theta\tau},\g^{\theta\tau}]$ admits a Cartan decomposition $\g'=\k'+\p'$ with $\k'\subset \k$, $\p'\subset \p$. Due to our choice, $\a^{-\tau}$ is a maximal abelian subalgebra of $\p'$. The set of roots $\Delta(\g',\a^{-\tau})$ consists exactly of of those reduced roots $\beta \in \Delta(\g,\a^{-\tau})$ for which $(\g^\beta)_+ \neq 0$. Moreover, the Weyl group ${\mathcal W'}$ associated to $\g'$ coincides with ${\mathcal W}_{K\cap H}$.
\end{remark}

We can now state the central theorem.

\begin{theorem}\label{VDB}{\em{(Van den Ban \cite{V})}} \\
Let $G$ be a real connected semisimple Lie group with finite center, equipped with an involution $\tau$, 
and $H$ a connected open subgroup of $G^\tau$. For $X \in \a^{-\tau}$, write $a=\exp X \in A^{-\tau}$. Then
\[ (pr_{\a^{-\tau}} \circ \log \circ \mu) ( H a ) = {\rm conv}(\mathcal{W}_{K \cap H}.X) - \Gamma(\Delta^{+}_{-}) . \]
\end{theorem}

\begin{remark} \ 
\begin{itemize}
\item The statement of the theorem above differs from the original in \cite{V} by a minus sign in front of the conal part $\Gamma(\Delta^{+}_{-})$. This is due to the fact that we consider the set $Ha$ and an Iwasawa decomposition $G=NAK$, whereas in \cite{V} the set $aH\subset G=KAN$ is considered. Indeed, if we denote the two middle projections by $\mu: NAK\to A$ and $\mu': KAN\to A$, then $\Gamma(\Delta^{+}_{-})=\log\circ \mu'(H) = - \log\circ\mu(H)$.
\item Van den Ban proved his theorem under the weaker condition that $H$ is an essentially connected 
open subgroup of $G^\tau$ (by reducing it to the connected case). 
\item If $\tau=\theta$ one obtains Kostant's (nonlinear) convexity theorem. Note that in this case 
the group $H$ and the orbit $Ha$ are compact.
\end{itemize}
\end{remark}

\section{Poisson structure}\label{Poisson}

Let $G$ be a connected and simply connected semisimple complex Lie group with Lie algebra $\g$. 
The Cartan involutions on both group and Lie algebra level will be denoted by $\theta$. In addition, 
let $\tau$ be a complex antilinear involution (on $G$ and $\g$) which commutes with $\theta$. 

The Lie algebra $\g$ decomposes into $(+1)$- and $(-1)$-eigenspaces with respect to both involutions 
$\theta$ and $\tau$.
\[ \g=\k+\p=\h+\q , \]
where $\k$ and $\h$ denote the $(+1)$-eigenspaces with respect to $\theta$ and $\tau$, respectively, 
and $\p$ and $\q$ denote the $(-1)$-eigenspaces. 

The maximal compact subgroup $K$ of $G$ with Lie algebra $\k$ is $\tau$-stable. Let $H$ denote the 
connected subgroup of $G$ consisting of $\tau$-fixed points. We will be interested in certain 
$H$-orbits in the symmetric space $G/K$. Each such orbit can be equipped with a Poisson structure 
as introduced by Evens and Lu. We briefly describe their method which can be found in \cite[Section 2.2]{EL}. 
For details on Poisson Lie groups see e.g. \cite{LW}. 

Let $(U,\pi_U)$ be a connected Poisson Lie group with tangent Lie bialgebra $(\uu,\uu^*)$ and double Lie algebra $\d=\uu \bowtie \uu^*$. The pairing
\[ \langle v_1+\lambda_1 , v_2+\lambda_2 \rangle := \lambda_1(v_2)+\lambda_2(v_1) \quad \forall \ v_1, v_2 \in \uu, \lambda_1, \lambda_2 \in \uu^* , \]
defines a non-degenerate symmetric bilinear form and turns $(\d,\uu,\uu^*)$ into a Manin triple. We will identify $\d^*$ with $\d$ via $\langle , \rangle $. 

Consider the following bivector $R \in \wedge^2 \d$:
\[ R(v_1+\lambda_1 , v_2+\lambda_2)=\lambda_2(v_1)-\lambda_1(v_2) \quad \forall \ v_1, v_2 \in \uu, \lambda_1, \lambda_2 \in \uu^* . \]
In terms of a basis $\{ v_1, \dots , v_n \}$ for $\uu$ and a dual basis $\{ \lambda_1, \dots , \lambda_n \}$ 
for $\uu^*$ the bivector is represented by $R=\sum^n_{i=1} \lambda_i \wedge v_i $. 

Assume that $D$ is a connected Lie group with Lie algebra $\d$, and assume that $U$ is a connected
subgroup of $D$ with Lie algebra $\uu$. Then 
$D$ acts on the Grassmannian ${\rm Gr}(n,\d)$ of $n$-dimensional subspaces of $\d$ via the adjoint action 
of $D$ on $\d$ and therefore defines a Lie algebra antihomomorphism
\[ \eta:\d \rightarrow \mathcal{X}({\rm Gr}(n,\d)) , \]
into the vector fields on ${\rm Gr}(n,\d)$. Using the symbol $\eta$ also for its multilinear extension we can define a bivector field $\Pi$ on ${\rm Gr}(n,\d)$ by
\[ \Pi = \frac{1}{2} \eta(R) . \]
Note that $\Pi$ in general does not define a Poisson structure on the entire ${\rm Gr}(n,\d)$. 
However, it does so on the subvariety $\mathfrak{L}(\d)$ of Lagrangian subspaces (with respect to $\langle , \rangle $ ) 
on $\d$, and on each $D$-orbit $D.\l \subset \mathfrak{L}(\d)$. 

The bivector $R$ also gives rise to a Poisson structure $\pi_-$ on $D$ that makes $(D,\pi_-)$ a Poisson Lie group:

\begin{equation}\label{DPoisson}
\pi_-(d)=\frac{1}{2} (r_d R-l_d R) \quad \forall d\in D . 
\end{equation}

Here $r_d$ and $l_d$ denote the differentials of right and left translations by $d$. 
Note that the restriction of $\pi_-$ to the subgroup $U\subset D$ coincides with the original Poisson structure 
$\pi_U$ on $U$, i.e. $(U,\pi_U)$ is a Poisson subgroup of $(D,-\pi_-)$. 

For $\l \in \mathfrak{L}(\d)$ the $D$-orbit through $l$ is not only a Poisson manifold with respect 
to $\Pi$ but a homogeneous Poisson space under the action of $(D,\pi_-)$. Moreover, the $U$-orbit 
$U.\l$ is a homogeneous $(U,\pi_U)$-space, since the Poisson tensor $\Pi$ at $\l$ turns out to be tangent to $U.\l$. 
In fact, the tangent space at $\l \in D.\l$ can be identified with $\d/ n(\l)$, where $n(\l)$ is the 
normalizer subalgebra of $\l$. In the case when $n(\l)=\l$, we identify the cotangent space with $\l$ itself, and for $X,Y \in \l$ one obtains:
\begin{equation}\label{distribution}
\Pi(\l)(X,Y)=\langle pr_\uu X,Y \rangle , \quad {\rm i.e.} \quad \Pi(\l)^\sharp(X)=pr_\uu X , 
\end{equation}
where $pr_\uu:\d \rightarrow \uu $ denotes the projection along $\uu^*$. 

Let $U^*$ be the connected subgroup of $D$ with Lie algebra $\uu^*$.  
What has been said about the Poisson Lie group $U$ is also true for its dual group $U^*$, 
i.e. $(U^*,\pi_{U^*})$ is a Poisson Lie subgroup of $(D,\pi_-)$ and the orbit $U^*.\l$ is a 
homogeneous $(U^*,\pi_{U^*})$-space. It follows in particular that $(U.\l)\cap (U^*.\l)$ 
contains the symplectic leaf through $\l$.

\vspace{.5cm}
\par
We now want to apply this construction to our complex semisimple Lie algebra $\g$. In the above notation 
we will have $\d=\g$, and the pairing $\langle , \rangle $ will be given by the imaginary part, $\Im \kappa$, of 
the Killing form $\kappa$ on $\g$. Note that $\k \in \mathfrak{L}(\d)$. Throughout the paper, we will identify the $G$-orbit 
through $\k$ with the symmetric space $G/K$. In particular, orbits in $G.\k$ are identified with those in $G/K$.  
Then we set $\uu=\h$, and it remains to define $\uu^*$. 

First we choose an appropriate Iwasawa decomposition of $\g$. Recall the $\tau$-stable Cartan decomposition $\g=\k+\p$. 
We fix a maximal abelian subalgebra $\a^{-\tau}$ in $\p \cap \q $. Then we can find an abelian subalgebra $\a^\tau$ 
in $\p \cap \h $ such that $\a=\a^{-\tau}+\a^\tau $ is maximal abelian in $\p$.
%The root space decomposition of $\g$ with respect to $\a$ is
%\[ \g=(\a+i\a)+\sum_{\alpha \in \Delta(\g,\a)} \g^\alpha . \]
%We choose a positive system $\Delta^+(\g,\a)$ of roots such that 
%$\alpha > \bet%a$ for all $\alpha, \beta \in \Delta^+(\g,\a)$ with 
%$\alpha |_{\a^{-\tau}}\neq %0$ and $\beta |_{\a^{-\tau}}=0$.
We choose a positive root system, $\Delta^+(\g, \a)$ by the lexicographic ordering with respect to an ordering of 
a basis of $\a$, which was constructed from a basis of $\a^{-\tau}$ followed by a basis of $\a^{\tau}$. 
This yields an Iwasawa decomposition $\g=\n+\a+\k$ which is compatible with the involution $\tau$ in the following sense.

\begin{lemma}\label{transversal}
For our choice of Iwasawa decomposition $\g=\n+\a+\k$, we have
\[ \h \cap \n = \{ 0 \} . \]
Besides, the centralizer of $\a^{-\tau}$ in $\g$ is a Cartan subalgebra of $\g$. 
\end{lemma}

\begin{proof}
Consider the root space decomposition of $\g$ with respect to $\a$,
\[ \g=(\a+i\a)+\sum_{\alpha \in \Delta(\g,\a)} \g^\alpha . \]
It is well-known \cite[Proposition 6.70]{Knapp} that there are no real roots for a maximally 
compact Cartan subalgebra $(i\a^{-\tau}+\a^{\tau})$ 
of $\h$, and therefore there are no $\alpha \in \Delta(\g,\a) $ such that $\alpha |_{\a^{-\tau}}= 0$. 
By \cite[Chapter VI, Lemma 3.3]{Helga}, this  implies 
that $\tau(\g^\alpha) \subset \bigoplus_{\alpha\in\Delta^+(\g, \a)}\g^{-\alpha}$ for 
all $\alpha \in \Delta^+(\g,\a)$, and the claim $\h \cap \n = \{ 0 \} $ follows immediately. 

Since each $\alpha\in \Delta(\g, \a)$ does not vanish outside a hyperplane of  
$\a^{-\tau}$, it follows that $\a^{-\tau}$ contains regular elements and its centralizer in $\g$
is a Cartan subalgebra of $\g$.  \\

\end{proof}

%Assume there is an $\alpha \in \Delta(\g,\a) $ such that $\alpha |_{\a^{-\tau}}= 0 $. 
%Then $\tau $ leaves the one dimensional complex vector space $\g^\alpha$ stable. 
%As an antilinear involution $\tau$ must have a $(+1)$-eigenvector in $\g^\alpha$. 
%Let $X_\alpha \in \g^\alpha$ be nonzero with $\tau(X_\alpha)=X_\alpha$. For $s\in \R$ 
%consider the group elements $u(s)=\exp(is(X_\alpha - \theta X_\alpha))\in Z_{K}(\a^{-\tau})$. 
%A direct calculation shows that for $H_\alpha =[X_\alpha, \theta X_\alpha ]$,
%$$ Ad(u(s)).H_\alpha = $$ $$ = 
%( \sum_{j=0}^{\infty} \frac{(2\alpha(H_\alpha))^j}{(2j)!}s^{2j}) H_\alpha - i \alpha(H_\alpha) 
%( \sum_{j=0}^{\infty} \frac{(2\alpha(H_\alpha))^j}{(2j+1)!}s^{2j+1}) (X_\alpha + \theta X_\alpha ). $$
%From $H_\alpha \in \p$ it follows that $\alpha(H_\alpha)\in \R$. 
%The non-degeneracy of the bilinear form $B(X,Y)=-\kappa(X,\theta Y)$ 
%implies $\alpha(H_\alpha)<0$. In particular, one can choose $X_\alpha$ 
%such that $\alpha(H_\alpha)=-1/2$. Setting $s=\pi /2$ yields
%\[ Ad(u(\pi /2)).H_\alpha = \frac{i}{2} (X_\alpha + \theta X_\alpha ). \]
%Note that $Ad(u(\pi /2)).H_\alpha \in \p \cap \z_{\g}(\a^{-\tau})$ and 
%$\tau(Ad(u(\pi /2)).H_\alpha)=-Ad(u(\pi /2)).H_\alpha$. 
%But this contradicts the maximality of $\a^{-\tau}$. 

Consider the Cartan subalgebra $\c=\z(\a^{-\tau})$ of $\g$. Lemma \ref{transversal} together with the properties of 
$\kappa$ implies that $\g=\h \oplus (\c^{-\tau}\oplus \n)$ is a Lagrangian splitting with respect to the 
bilinear form $\Im \kappa$. In other words, $\ (\g,\h,(\c^{-\tau}+\n)) \ $ is a Manin triple. 

We can now define the desired Poisson manifolds using the method of Evens and Lu outlined above. We set
\[ \d=\g, \ \uu=\h, \ \uu^*=\c^{-\tau}+\n, \ \langle , \rangle =\Im \kappa . \]
%Note that $\k\in \mathfrak{L}(\g)$ and $G.\k\cong G/K$. 
Let $C$, $C^{-\tau}$, $A$ and $N$ denote the analytic subgroups 
of $G$ with Lie algebras $\c$, $\c^{-\tau}$, $\a$ and $\n$, respectively. 
The group $H$ now has the structure of a Poisson Lie group. Its dual group is $H^*=C^{-\tau}N$.
Fix $a\in A^{-\tau}$ and consider the base point $a.K\in G/K$.  
The $H$-orbit $P_a=Ha.K\in G/K$ is a Poisson homogeneous manifold with respect to 
the action by $(H,\pi_H)$. Also, the dual group orbit $H^*a.K$ is Poisson homogeneous with respect to $\pi_{H^*}$. 
For the symplectic leaf in $P_a$ through $a$, denoted by $M_a$, we have 
$M_a \subseteq Ha.K\cap H^*a.K$.

\begin{lemma}\label{fibration}
The Poisson manifold $P_a$ is regular and equals the union of  
$A^{\tau}$-translates of $M_a$, i.e. each $p\in P_a$ can be 
written $p=a'm$ with unique $a'\in A^{\tau}, m\in M_a$. Moreover, $M_a = Ha.K\cap H^*a.K$.
\end{lemma}

\begin{proof}
Consider the map $M:A^\tau \times M_a \rightarrow P_a $. 

First we will show that $M$ is injective.
The Poisson tensor $\pi_H=\pi_-$ as defined in (\ref{DPoisson}) vanishes at each element 
$c\in C^\tau$, since $Ad(c)$ leaves both $\h$ and $\h^*=\c^{-\tau}+\n$ stable. 
Therefore $a'\in A^\tau $ acts on $P_a$ by Poisson diffeomorphisms and 
maps the symplectic leaf $M_a$ onto the symplectic leaf $M_{a'a}$. 
But $M_{a_1a}\neq M_{a_2a}$ for $a_1 \neq a_2 \in A^{\tau}$, 
following from the fact that $M_{a_1a}$ lies in $H^*a_1a.K=C^{-\tau}Na_1a.K$ and the 
uniqueness of the Iwasawa decomposition. 

At each point $p\in P_a$ one can explicitly calculate the codimension of the symplectic 
leaf through $p$ in $P_a$, for instance by means of an infinitesimal version of 
Corollary 7.3 in \cite{Lu} and Theorem 2.21 in \cite{EL}. It follows that the codimension
of the leaf through the point $p=ha.K$ in the orbit $P_a$ equals the dimension of the intersection of 
$Ad(a)\k$ and $Ad(h^{-1})\h^*$, which is easily seen to be  
independent from the point $p\in P_a$ and equal to the dimension of $\a^\tau$.
Here we used the fact that the dimension of $Ad(ha)\k\cap\h^*$ cannot exceed the dimension 
of $\a^{\tau}$, since the Killing form is negative definite on $Ad(ha)\k$ and a maximal negative 
definite subspace of $\h^*$ is $i\a^{\tau}$.  
This shows that $P_a$ is a regular Poisson manifold, and that $A^{\tau}M_a$ is a full dimensional 
subset of $P_a$. Since $A^{\tau}$ acts freely on $P_a$ and $P_a$ is a regular Poisson manifold, it 
can be represented as the union of such open subsets.  The connectedness of $P_a$ then implies 
that $P_a=A^\tau M_a$. 

Since $A^{\tau}$ is connected and the union of $A^{\tau}$-translates of $Ha.K\cap H^*a.K$ equals 
$Ha.K$ and thus is also connected, it is easy to see that $Ha.K\cap H^*a.K$ is connected as well. 
Besides, from the transversality we see that
$$
{\rm dim}(Ha.K\cap H^*a.K)={\rm dim}(Ha.K)+{\rm dim}(H^*a.K)-{\rm dim}(G/K).
$$
Note that the first part of the proof implies that $A^\tau a.K\cap H^*a.K =\{ a.K \} $.
Therefore, the codimension of $Ha.K\cap H^*a.K$ in $Ha.K$ is at least ${\rm dim}(\a^\tau)$. 
But since $M_a$ has codimension equal to ${\rm dim}(\a^\tau)$, 
and $M_a \subseteq Ha.K\cap H^*a.K$, the last inclusion is actually an equality.

\end{proof}

Consider the torus $T=\exp(i\a^{-\tau})\subset H$. It acts on $M_a$ in a symplectic manner, 
since $\pi_H$ vanishes at each $t\in T$. Moreover, the next lemma shows that this action 
is Hamiltonian with an associated moment map that is closely related to the middle 
projection $\mu:G=NAK \rightarrow A$ of the Iwasawa decomposition.

\begin{lemma}\label{mm}
The action of $\ T=\exp(i\a^{-\tau})$ on $M_a$ is Hamiltonian with a moment map 
$\Phi=pr_{\a^{-\tau}} \circ \log \circ \mu$.
Here, $pr_{\a^{-\tau}}:\a \rightarrow \a^{-\tau}$ denotes the projection along $\a^\tau$, and $\t^*$
is identified with $\a^{-\tau}$ via $\Im \kappa$. \\
Moreover, the moment map $\Phi$ is proper.
\end{lemma}
\begin{proof}
\begin{enumerate}
\item $\Phi=pr_{\a^{-\tau}}\circ \log \circ \mu$ is a moment map. 

Let $b:G=NAK \longrightarrow B=NA$ be the $B$-projection in the Iwasawa decomposition. 
We write $pr_\a:\g=\n+\a+\k \rightarrow \a $ for the middle projection on the Lie algebra level. 
Let $Z\in \t=i\a^{-\tau}$, $h\in  H$ and $X\in \h$. We denote by $\Phi_Z$ the function obtained by 
evaluating $\Phi$ at $Z$, by ${\tilde X}_{ha}$ the tangent vector of the vector field generated by $X$ 
at the point $ha.K\in M_a$ (for brevity we will write $h.K$ simply as $h$ henceforth, 
without fear of confusion) and by $D\Phi_{b(ha)}$ the derivative of $\Phi$ at the point $b(ha)$. 
%regard $M_a$ as a subset of $G/K$, however we can treat $\Phi$ also as a map 
%from $M_a$ to $\a^{-\tau}$, because of the natural identification of $NA$ with $G/K$. 
We have:

\begin{eqnarray*}
d\Phi_Z(ha).\tilde{X}_{ha} & = & \frac{d}{ds}\Big|_{s=0} \Phi_Z(\exp(sX)ha) = \langle \frac{d}{ds}\Big|_{s=0} \Phi(\exp(sX)ha),Z \rangle \\
& = & \langle \frac{d}{ds}\Big|_{s=0} \Phi(b(ha)\exp(sAd(b(ha)^{-1})X)),Z \rangle \\ & = &
\langle D\Phi_{b(ha)}Ad(b(ha)^{-1})X,Z \rangle \\
& = & \langle pr_{\a^{-\tau}} \circ pr_\a Ad(b(ha)^{-1})X,Z \rangle = \langle Ad(b(ha)^{-1})X,Z \rangle \\
& = & \langle X, Ad(b(ha))Z \rangle
\end{eqnarray*}

The second last step follows from the fact that $\t$ and $\k+\a^\tau+\n$ are orthogonal with respect to $\langle , \rangle $. 

Note that $Ad(b(ha))Z \in Z+\n$. With (\ref{distribution}) this implies
\[ \Pi(ha)^{\sharp}(d\Phi_Z(ha)) = pr_{\h}Ad(b(ha))Z = Z . \]

\item $\Phi$ is proper. 

This follows from Lemma 3.3 in \cite{V}, which states the properness of the map
\[ F_a:(H\cap L_0)\backslash H \rightarrow \a^{-\tau}, \quad F_a(x)=\Phi(xa) . \]
In our case $L_0=\exp(i\a )A^\tau$ (since $\z_{\g}(\a^{-\tau})=\c$ by the argument in the proof of Lemma \ref{transversal}). \\
Properness of the map $F_a: TA^{\tau} \backslash H \rightarrow \a^{-\tau}$ implies properness of the induced maps 
$F_a:A^{\tau} \backslash H \rightarrow \a^{-\tau}$ and $F_a: A^{\tau} \backslash H / (H\cap aKa^{-1}) \rightarrow \a^{-\tau}$. 
Since $A^{\tau} \backslash H / (H\cap aKa^{-1}) \cong M_a$ by Lemma \ref{fibration}, and since $F_a$ becomes $\Phi$ 
under this identification, the claim follows.
\end{enumerate}
\end{proof}

\begin{remark}
In case $\tau=\theta$ the Lu-Evens Poisson structure on $P_a=Ka.K$ coincides with the
Lu-Weinstein symplectic structure, and Lemma \ref{mm} becomes Theorem 4.13 in \cite{LR}.
\end{remark}

\section{Symplectic convexity}\label{sympl}

Throughout this section we assume $G$ to be complex and the involution $\tau$ to be complex antilinear. In this case we will interpret van den Ban's theorem in the symplectic framework developed in Section \ref{Poisson}. More precisely, it can be viewed as a corollary of a symplectic convexity theorem for Hamiltonian torus actions.

Van den Ban's theorem describes the image of the group orbit $Ha$ under the map $ pr_{\a^{-\tau}} \circ \log \circ \mu $. Recall from Section \ref{Poisson} the symplectic manifold $M_a \subseteq Ha.K \subseteq G/K$  on which the torus $T=\exp(i\a^{-\tau})$ acts in a Hamiltonian fashion (Lemma \ref{mm}). The associated moment map is $\Phi=pr_{\a^{-\tau}} \circ \log \circ \mu $. From Lemma \ref{fibration} and from the $A^\tau$-invariance of $ pr_{\a^{-\tau}} \circ \log \circ \mu $ it follows that
\[ (pr_{\a^{-\tau}} \circ \log \circ \mu)(Ha)=\Phi(M_a). \]
This means that van den Ban's theorem can be viewed as a description of the image of a symplectic manifold under an appropriate moment map.
   
The description of the image of the moment map is the content of a series of symplectic convexity theorems. Probably best known are the original theorems of Atiyah and Guillemin-Sternberg \cite{A,GS}. The result needed here is a generalization of the AGS-theorems to a non-compact setting. Several versions can be found in the literature, e.g. \cite{LMTW,Pr}. We will state the theorem as given in \cite{HNP}. Recall that a subset $C$ of a finite dimensional vector space $V$ is called locally polyhedral iff for each $x\in C$ there is a neighborhood $U_x\subseteq V$ such that $C \cap U_x = (x+\Gamma_x) \cap U_x$ for some cone $\Gamma_x$. A cone $\Gamma$ is called proper if it contains no lines, otherwise $\Gamma$ is called improper.

\begin{theorem}\label{hnp} {\em{\cite[Theorem 4.1(i)]{HNP}}} 
Consider a Hamiltonian torus action of $T$ on the connected symplectic manifold $M$. Suppose the associated moment map $\Phi:M \rightarrow \t^*$ is proper, i.e. $\Phi$ is a closed mapping and $\Phi^{-1}(Z)$ is compact for every $Z\in \t^*$. 
Then $\Phi(M)$ is a closed, locally polyhedral, convex set.
\end{theorem}

\begin{remark}\label{hnpremark}
Theorem 4.1 in \cite{HNP} contains more detailed information, in particular a description of the cones that span $\Phi(M)$ locally (part (v)). More precisely, for each $m\in M$ there is a neighborhood $U_{\Phi(m)}\subseteq \t^*$ of $\Phi(m)$ such that $\Phi(M) \cap U_{\Phi(m)} = (\Phi(m)+\Gamma_{\Phi(m)}) \cap U_\Phi(m)$, where $\Gamma_{\Phi(m)}=\t_m^\perp + C_m$. Here, $\t_m$ denotes the Lie algebra of the stabilizer $T_m$ of $m$, and $C_m\subseteq \t_m^*$ denotes the cone which is spanned by the weights of the linearized action of $T_m$. The (nontrivial) fact that the cone $\Gamma_{\Phi(m)}=\t_m^\perp + C_m$ is actually independent of the choice of a preimage point of $\Phi(m)$ is also shown in \cite{HNP}.
\end{remark}

Coming back to the symplectic manifold $M_a$, Lemma \ref{mm} shows that the moment map $\Phi=pr_{\a^{-\tau}} \circ \log \circ \mu $ on $M_a$ is proper. Theorem \ref{hnp} can therefore be applied and yields
\[ \Phi(M_a) \quad \mbox{is a closed, locally polyhedral, convex set.} \]

We will now give a more detailed description of $\Phi(M_a)$. It turns out that the $T$-action on $M_a$ has (finitely many) fixed points. At each fixed point we can calculate the cones that locally span $\Phi(M_a)$. From this description it will follow that the entire set $\Phi(M_a)$ lies in a proper cone and can therefore be described entirely by the local data at the fixed points.

We begin by determining the $T$-fixed points.

\begin{prop}\label{fixed}
The $T$-fixed points in $M_a$ are exactly those elements of the form $w(a).K\in G/K$ with $w\in {\mathcal W}_{K\cap H}=N_{K \cap H}(\a^{-\tau})/Z_{K \cap H}(\a^{-\tau})$.
\end{prop}
\begin{proof}
Recall that for $a\in A^{-\tau}$ we view the symplectic manifold $M_a$ as a submanifold of the $H$-orbit in $G/K$ through the base point $a.K \in G/K$. Clearly, each element $w(a).K \in G/K$ with $w\in {\mathcal W}_{K\cap H}$ is $T$-fixed. To see that $w(a).K$ lies in $M_a$, note that $w(a).K\in H^*a.K$ since $w(a)\in A^{-\tau}$. On the other hand, there exists $k\in K\cap H$ such that $w(a)=kak^{-1}$, which implies $w(a).K\in Ha.K$. Therefore, $w(a).K\in Ha.K\cap H^*a.K =M_a$ by Lemma \ref{fibration}.

Conversely, assume that $cpa.K \in M_a$ with $c\in K^\tau , p\in \exp(\p^\tau)$ is $T$-fixed. Since $M_a$ lies in the orbit of the dual group $H^*=NC^{-\tau}$ there are elements $n\in N , b\in A^{-\tau}, k\in K$ such that $cpa=nbk$. Since $nb.K\in G/K$ is a $T$-fixed point,
\[ tnt^{-1}b \in nbK \qquad \forall \ t\in T . \]
The Lie subalgebra $\n$ is $T$-invariant, so by the uniqueness of the Iwasawa decomposition $tnt^{-1}=n$ for all $t\in T$. But since $\alpha |_{\a^{-\tau}} \neq 0$ for all $\alpha \in \Delta(\g,\a)$ this can happen only for $n=e$. This implies $cpa=bk$. 

Symmetrizing the last equation yields
\begin{equation}\label{fix1}
cpa \theta(cpa)^{-1}=cpa^2pc^{-1}=b^2 .
\end{equation}
Applying $\theta \circ \tau$ to (\ref{fix1}) gives
\begin{equation}\label{fix2}
cp^{-1}a^2p^{-1}c^{-1}=b^2 .
\end{equation}
We multiply (\ref{fix1}) by (\ref{fix2}) from the right and from the left and obtain
%\begin{equation}\label{fix3}
\[ cpa^4p^{-1}c^{-1}=b^4=cp^{-1}a^4pc^{-1} . \]
%\end{equation}
But then $pa^4p^{-1}=p^{-1}a^4p$, i.e. $p^2$ and $a^4$ commute (and are self-adjoint). Therefore, $p$ and $a^2$ also commute, and we can combine equations (\ref{fix1}) and (\ref{fix2}) to
\[ cp^2a^2c^{-1}=b^2=cp^{-2}a^2c^{-1}. \]
This shows $p^2=p^{-2}$ or $p=e$. 
But then (\ref{fix2}) implies $cac^{-1}=b$. Since both $a$ and $b$ lie in $A^{-\tau}$ and since $c\in K^\tau=K\cap H$, there is some element $w\in {\mathcal W}_{K\cap H}$ such that $w(a)=b$ (Recall from Remark \ref{thetatau} that ${\mathcal W}_{K\cap H}$ is the Weyl group of the reductive Lie algebra $\g^{\theta \tau}=(\k\cap \h)+(\p \cap \q)$ of $\theta \tau$-fixed points of $\g$).

The $T$-fixed point $cpa.K \in M_a$ can therefore be written as $cpa.K=b.K=w(a).K$ for some $w\in {\mathcal W}_{K\cap H}$. 
\end{proof} 

Recall our choice of base point $a=\exp(X)$ and the identification $\t^* \cong \a^{-\tau}$. We now describe the image of the moment map $\Phi(M_a)\in \a^{-\tau}$ in the neighborhood of a fixed point image $\Phi(w(a).K)=w(X)$. From Theorem \ref{hnp} (and Remark \ref{hnpremark}) we know that locally $\Phi(M_a)$ looks like $w(X)+\Gamma_{w(X)}$ for some cone $\Gamma_{w(X)} \in \a^{-\tau}$. The next Lemma describes $\Gamma_{w(X)}$ in terms of the vectors $H_\beta$ for (reduced) roots $\beta\in \Delta(\g,\a^{-\tau})$ defined in Section \ref{vdb}. 

\begin{lemma}\label{cones}
Let $a=\exp X$ with $X\in \a^{-\tau}$ and $w\in {\mathcal W}_{K\cap H}$. The local cone $\Gamma_{\Phi(w(a).K)}=\Gamma_{w(X)} \subseteq \a^{-\tau}$ is the cone spanned by the union of the following two sets.
\begin{eqnarray*}
\{ -\beta(w(X)) H_\beta : \beta \in \Delta^+(\g,\a^{-\tau}), (\g^\beta)_+ \neq 0 \} \\
{\rm and} \quad \{ -H_{\beta}: \beta \in \Delta^+(\g,\a^{-\tau}), (\g^\beta)_- \neq 0 \} 
\end{eqnarray*}
\end{lemma}

\begin{proof}
We are adapting the argument from \cite[page 155]{HNP} to our setting. 
To determine the local cone $\Gamma_{w(X)}$ it is enough to consider the linearized 
action of $T$ on the tangent space $V_{w(a).K}:=T_{w(a).K}M_a$. Darboux's theorem 
guarantees the existence of a $T$-equivariant symplectomorphism of a 
neighborhood of $w(a).K\in M_a$ onto a neighborhood of $0\in V_{w(a).K}$. This leads 
to a local normal form for the moment map.
\[ \Phi_Z(Y)=\frac{1}{2}\Omega_{w(a).K}(Z.Y,Y) \quad \forall \ Y\in V_{w(a).K}, Z\in \t . \]
Here, $\Omega_{w(a).K}$ denotes the symplectic form on the symplectic vector space $V_{w(a).K}$. Since $T$ acts symplectically on $V_{w(a).K}$ the notation $Z.Y$ makes sense as the linear action of an element $Z\in \mathfrak{sp}(V_{w(a).K})$ on a vector $Y\in V_{w(a).K}$.
In appropriate symplectic coordinates $q_1, p_1,\dots,q_n, p_n$ we have $\Omega_{w(a).K}=\sum_i dq_i\wedge dp_i$ and the matrix representation for the linear map defined by $Z\in \t$ is
\[ Z.(q_1,p_1,\dots,q_n,p_n)=\begin{pmatrix} 0 & \alpha_1(Z) & & & \\ -\alpha(Z) & 0 & & & \\ & & \ddots & & \\ & & & 0 & \alpha_n(Z) \\ & & & -\alpha_n(Z) & 0 \end{pmatrix} \begin{pmatrix} q_1 \\ p_1 \\ \vdots \\ q_n \\ p_n \end{pmatrix} . \]
The moment map takes the form
\[ \Phi(q_1,p_1,\dots,q_n,p_n)=\Phi(w(a).K)+\sum_{i=1}^n \alpha_i \frac{1}{2}(q_i^2+p_i^2) . \]
In terms of the symplectic coordinates on $V_{w(a).K}$ chosen above, the matrix representations for the symplectic form $\Omega_{w(a).K}$ and the corresponding Poisson tensor $\Pi_{w(a).K}$ just differ by a factor of $(-1)$. The moment map can then be expressed in terms of the Poisson tensor.
\[ \Phi_Z(\varphi)=-\Pi_{w(a).K}(Z.\varphi,\varphi) \qquad \forall \ \varphi \in V_{w(a).K}^*, Z\in \t . \]
(Recall the bijection $\Pi^\sharp:V_{w(a).K}^* \rightarrow V_{w(a).K}$. Then $Z.\varphi=(\Pi^\sharp)^{-1}(Z.(\Pi^\sharp(\varphi)))$, where the dot on the right hand side has been explained above.)

The local cone $\Gamma_{w(X)}$ is just $\Phi(V_{w(a).K}^*)$, i.e. it consists exactly of the weights
\begin{equation}\label{weight}
\{ \ Z\mapsto -\Pi_{w(a).K}(Z.\varphi,\varphi) \ : \ \varphi \in V_{w(a).K}^* \ \} 
\end{equation}
Recall that we identify the cotangent space $T_{w(a).K}^*(G/K)$ with $Ad(w(a)).\k $. The formula for the Poisson tensor at $w(a).K$ says that for $Y_1, Y_2 \in \k$,
\[ \Pi_{w(a).K}(Ad(w(a))Y_1,Ad(w(a))Y_2) = \langle pr_\h Ad(w(a))Y_1,Ad(w(a))Y_2 \rangle . \]
Note that $T_{w(a).K}^*(G/K)=T_{w(a).K}^*M_a \oplus (T_{w(a).K}M_a)^\perp$. Both $T_{w(a).K}^*M_a$ and $(T_{w(a).K}M_a)^\perp$ are stable under the action of $T$. Moreover, $T_{w(a).K}M_a=\Pi_{w(a).K}^{\sharp}(T_{w(a).K}^*(G/K))$ by the definition of the symplectic leaf $M_a$. Hence, for $\varphi \in T_{w(a).K}^*M_a, \psi \in (T_{a_w}M_{a_w})^\perp$ and $Z\in \t $, one obtains
\begin{eqnarray*}
\Pi_{w(a).K}(Z.(\varphi+\psi),(\varphi+\psi)) & = & (\varphi+\psi).\Pi_{w(a).K}^{\sharp}(Z.(\varphi+\psi)) \\
& = & \varphi.\Pi_{w(a).K}^{\sharp}(Z(\varphi+\psi)) \\
& = & \Pi_{w(a).K}(Z\varphi+Z\psi,\varphi) \\ & = & -(Z\varphi+Z\psi).\Pi_{w(a).K}^{\sharp}(\varphi) \\
& = & \Pi_{w(a).K}(Z\varphi,\varphi) 
\end{eqnarray*}

In view of (\ref{weight}) and (\ref{distribution}) (from Section \ref{Poisson}) it follows that the local cone is given by
\begin{equation}\label{cone}
\Gamma_{w(X)} = \{ \ Z\mapsto -\langle pr_\h [Z,Ad(w(a))Y],Ad(w(a))Y \rangle \ : \ Y\in \k \ \} .
\end{equation}
In order to determine the weights in (\ref{cone}) we will construct a basis $\{ v_1,\dots,v_r\}$ for $\k$ with two main features.
\begin{enumerate} 
\item For each $v_i$ we determine explicitly an element $H_i\in \a^{-\tau}$ such that 
\[ \langle pr_\h [Z,Ad(w(a))v_i],Ad(w(a))v_i \rangle = \Im \kappa(H_i,Z) \quad \forall \ Z\in \t . \]
\item $\langle pr_\h [Z,Ad(w(a))v_i],Ad(w(a))v_j \rangle = 0$ for all $Z\in \t $ whenever $i\neq j$.
\end{enumerate}
Once such a basis is found each $Y\in \k$ can be written as a linear combination $Y=\sum_{i=1}^N c_i v_i$. Then, for $Z\in \t $,
\begin{eqnarray*}
\langle pr_\h [Z,Ad(w(a))Y],Ad(w(a))Y \rangle & = & \langle pr_\h [Z,Ad(w(a))\sum_{i=1}^N c_i v_i],Ad(w(a))\sum_{i=1}^N c_i v_i \rangle \\
& = & \sum_{i=1}^N c_i^2 \langle pr_\h [Z,Ad(w(a))v_i],Ad(w(a))v_i \rangle \\
& = & \sum_{i=1}^N c_i^2 \Im \kappa(H_i,Z)
\end{eqnarray*}
In view of (\ref{cone}) it then follows that $\Gamma_{w(X)}$ is the cone spanned by the vectors $H_i$.

Recall the weight space decomposition of $\g$ with respect to $\a^{-\tau}$.
\[ \g=\a^{-\tau} \oplus \a^\tau \oplus i\a^{-\tau} \oplus i\a^\tau \oplus \sum_{\beta \in \Delta(\g,\a^{-\tau})} \g^\beta \]
Each $\g^\beta$ is stable under the involution $\theta \tau$, hence decomposes into $(+1)$- and $(-1)$-eigenspaces $\g^\beta=(\g^\beta)_+ \oplus (\g^\beta)_-$. We first consider certain bases for $\g^\beta=(\g^\beta)_+$ and $\g^\beta=(\g^\beta)_-$. Each $\g^\beta$ is stable under the adjoint action of $\a^\tau$. For the corresponding weight space decomposition we write
\[ \g^\beta = \sum_{\eta \in \Delta(\g^\beta,\a^\tau)} \g^{\beta,\eta} \]
Note that $\g^{\beta,\eta}$ is equal to the eigenspace $\g^\alpha \subset \n$ for $\alpha \in \Delta(\g,\a)$ if and only if $\alpha |_{\a^{-\tau}}=\beta$ and $\alpha |_{\a^\tau}=\eta$. Also, if $\g^{\beta,\eta}=\g^\alpha$, then $\beta \in \Delta^+(\g,\a^{-\tau})$ if and only if $\alpha \in \Delta^+(\g,\a)$. The involutions $\tau$ and $\theta$ transform the eigenspaces as follows
\[ \tau(\g^{\beta,\eta})=\g^{-\beta,\eta}, \quad \theta(\g^{\beta,\eta})=\g^{-\beta,-\eta}, \quad \theta\tau(\g^{\beta,\eta})=\g^{\beta,-\eta} \]
For each eigenspace $\g^{\beta,\eta}$ fix a vector $X_{\beta,\eta}$ that spans $\g^{\beta,\eta}$ as a complex vector space. If $\eta\neq 0$ we define
\[ A_{\beta,\eta}=X_{\beta,\eta}+\theta\tau X_{\beta,\eta}, \quad B_{\beta,\eta}=X_{\beta,\eta}-\theta\tau X_{\beta,\eta} . \]
We obtain the following (complex) basis for the reduced root space $\g^\beta$ 
\[ \{ X_{\beta,0} \} \cup \{ A_{\beta,\eta}:\eta \neq 0 \} \cup \{ B_{\beta,\eta}:\eta \neq 0 \} \]
The important feature of this basis is that it consists of eigenvectors of the complex linear involution $\theta\tau$. Indeed, $\theta\tau A_{\beta,\eta}=A_{\beta,\eta}, \theta\tau B_{\beta,\eta}=-B_{\beta,\eta}$ and $X_{\beta,0}$ might be a $(+1)$- or a $(-1)$-eigenvector of $\theta\tau$. Therefore, a basis for $(\g^\beta)_+$ is given by the $A_{\beta,\eta}$'s and possibly $X_{\beta,0}$. A basis for $(\g^\beta)_-$ is given by the $B_{\beta,\eta}$'s and possibly $X_{\beta,0}$ (iff it is not contained in $\g^\beta=(\g^\beta)_+$).

The desired (real) basis for $\k$ now consists of a basis for $\z_\k(\a)=\z_\k(\a^{-\tau})=i\a^{-\tau}+i\a^\tau $ and the following set.
\begin{equation}\label{kbase}
\begin{split}
\bigcup_{\beta \in \Delta^+(\g,\a^{-\tau})} & \big( \{ X_{\beta,0}+\theta X_{\beta,0} \} \cup \{ iX_{\beta,0}+\theta iX_{\beta,0} \} \\
& \ \cup \{ A_{\beta,\eta}+\theta A_{\beta,\eta}: \eta \neq 0 \} \cup \{ iA_{\beta,\eta}+\theta iA_{\beta,\eta}: \eta \neq 0 \} \\
& \ \cup \{ B_{\beta,\eta}+\theta B_{\beta,\eta}: \eta \neq 0 \} \cup \{ iB_{\beta,\eta}+\theta iB_{\beta,\eta}: \eta \neq 0 \} \ \big)
\end{split}
\end{equation}
We can now calculate the weights appearing in (\ref{cone}) for each basis element. We fix $Z=iH \in \t=i\a_{-\tau}$. Recall that $a=\exp X$, therefore $w(a)=\exp(w(X))$. First we make two short auxiliary calculations. For a vector $C_\beta \in \g^\beta$ which is also a $\theta\tau$-fixed point,
\begin{eqnarray*}
[Z,Ad(w(a)).(C_\beta +\theta C_\beta)] & = & i\beta(H) w(a)^\beta C_\beta - i\beta(H) w(a)^{-\beta} \theta C_\beta \\
& = & \beta(H) w(a)^{-\beta} (iC_\beta+\theta i C_\beta)+\beta(H)(w(a)^\beta-w(a)^{-\beta}) i C_\beta 
\end{eqnarray*}
In the second line, the first summand lies in $\h$ the second in $\c^{-\tau}+\n$. For $D_\beta \in \g^\beta$ such that $\theta\tau D_\beta=-D_\beta$, the $\h\oplus (\c^{-\tau}+\n)$ decomposition is different:
\begin{eqnarray*}
[Z,Ad(w(a)).(D_\beta +\theta D_\beta)] & = & i\beta(H) w(a)^\beta D_\beta - i\beta(H) w(a)^{-\beta} \theta D_\beta \\
& = & \beta(H) w(a)^{-\beta} (-iD_\beta+\theta i D_\beta)+\beta(H)(w(a)^\beta+w(a)^{-\beta}) i D_\beta 
\end{eqnarray*}

Now, for $A_{\beta,\eta}$, which lies in $\g^\beta$ and satisfies $\theta\tau A_{\beta,\eta}=A_{\beta,\eta}$, we compute 
\begin{multline}\label{Aweight}
\langle pr_\h [Z,Ad(w(a)).(A_{\beta,\eta}+\theta A_{\beta,\eta})],Ad(w(a)).(A_{\beta,\eta}+\theta A_{\beta,\eta}) \rangle \\
\begin{split}
= & \ \langle \beta(H) w(a)^{-\beta} (iA_{\beta,\eta}+\theta i A_{\beta,\eta}),w(a)^\beta A_{\beta,\eta}+w(a)^{-\beta}\theta A_{\beta,\eta} \rangle \\
= & \ \beta(H) w(a)^{-2\beta} \langle iA_{\beta,\eta},\theta A_{\beta,\eta} \rangle + \beta(H) \langle \theta i A_{\beta,\eta},A_{\beta,\eta} \rangle \\
= & \ (w(a)^{-2\beta}-1) \ \Re \kappa (A_{\beta,\eta},\theta A_{\beta,\eta}) \ \beta(H) \\
= & \ (w(a)^{-2\beta}-1) \ \Re \kappa (A_{\beta,\eta},\theta A_{\beta,\eta}) \ \kappa(H_\beta,H) \\
= & \ (w(a)^{-2\beta}-1) \ \Re \kappa (A_{\beta,\eta},\theta A_{\beta,\eta}) \ \Im \kappa(H_\beta,Z)
\end{split}
\end{multline}
We can replace $A_{\beta,\eta}$ with $iA_{\beta,\eta}$ in the above calculation and obtain
\begin{multline*}
\langle pr_\h [Z,Ad(w(a)).(iA_{\beta,\eta}+\theta iA_{\beta,\eta})],Ad(w(a)).(iA_{\beta,\eta}+\theta iA_{\beta,\eta}) \rangle \\
\begin{split}
= & \ (w(a)^{-2\beta}-1) \ \Re \kappa (iA_{\beta,\eta},\theta iA_{\beta,\eta}) \ \beta(H) \\
= & \ (w(a)^{-2\beta}-1) \ \Re \kappa (A_{\beta,\eta},\theta A_{\beta,\eta}) \ \Im \kappa(H_\beta,Z)
\end{split}
\end{multline*}
Carrying out the calculation for $B_{\beta,\eta}$ (which are $(-1)$-eigenvectors of $\theta\tau$) we obtain a result of a different nature
\begin{multline}\label{Bweight} 
\langle pr_\h [Z,Ad(w(a)).(B_{\beta,\eta}+\theta B_{\beta,\eta})],Ad(w(a)).(B_{\beta,\eta}+\theta B_{\beta,\eta}) \rangle \\
= -(w(a)^{-2\beta}+1) \ \Re \kappa (B_{\beta,\eta},\theta B_{\beta,\eta}) \ \Im \kappa(H_\beta,Z) , 
\end{multline}
and
\begin{multline*} 
\langle pr_\h [Z,Ad(w(a)).(iB_{\beta,\eta}+\theta iB_{\beta,\eta})],Ad(w(a)).(iB_{\beta,\eta}+\theta iB_{\beta,\eta}) \rangle \\
= -(w(a)^{-2\beta}+1) \ \Re \kappa (B_{\beta,\eta},\theta B_{\beta,\eta}) \ \Im \kappa(H_\beta,Z) . 
\end{multline*}
If $X_{\beta,0}$ is fixed by $\theta\tau$, then
\begin{multline}\label{X+weight} 
\langle pr_\h [Z,Ad(w(a)).(X_{\beta,0}+\theta X_{\beta,0})],Ad(w(a)).(X_{\beta,0}+\theta X_{\beta,0}) \rangle \\
= (w(a)^{-2\beta}-1) \ \Re \kappa (X_{\beta,0},\theta X_{\beta,0}) \ \Im \kappa(H_\beta,Z)  , 
\end{multline}
and
\begin{multline*} 
\langle pr_\h [Z,Ad(w(a)).(iX_{\beta,0}+\theta iX_{\beta,0})],Ad(w(a)).(iX_{\beta,0}+\theta iX_{\beta,0}) \rangle \\
= (w(a)^{-2\beta}-1) \ \Re \kappa (X_{\beta,0},\theta X_{\beta,0}) \ \Im \kappa(H_\beta,Z) . 
\end{multline*}
The case that $\theta\tau X_{\beta,0}=-X_{\beta,0}$ leads to
\begin{multline}\label{X-weight} 
\langle pr_\h [Z,Ad(w(a)).(X_{\beta,0}+\theta X_{\beta,0})],Ad(w(a)).(X_{\beta,0}+\theta X_{\beta,0}) \rangle \\
= -(w(a)^{-2\beta}+1) \ \Re \kappa (X_{\beta,0},\theta X_{\beta,0}) \ \Im \kappa(H_\beta,Z)  ,\end{multline}
and
\begin{multline*} 
\langle pr_\h [Z,Ad(w(a)).(iX_{\beta,0}+\theta iX_{\beta,0})],Ad(w(a)).(iX_{\beta,0}+\theta iX_{\beta,0}) \rangle \\
= -(w(a)^{-2\beta})+1) \ \Re \kappa (X_{\beta,0},\theta X_{\beta,0}) \ \Im \kappa(H_\beta,Z) . \end{multline*}
Moreover, for $Y\in \z_\k(\a)$ one easily checks that
\[ \langle pr_\h [Z,Ad(w(a)).Y],Ad(w(a)).Y \rangle = 0 . \]

Note that the coefficient of $\Im \kappa(H_\beta,Z)$ in (\ref{Bweight}) and (\ref{X-weight}) is always positive. Therefore, basis vectors of $\k$ which are $(-1)$-eigenvectors of $\theta\tau$ contribute the set $\{ -H_{\beta}: \beta \in \Delta^+(\g,\a^{-\tau}), (\g^\beta)_- \neq 0 \}$ to $\Gamma_{w(X)}$.

On the other hand, the coefficient of $\Im \kappa(H_\beta,Z)$ in (\ref{Aweight}) and (\ref{X+weight}) depends on the value of $\beta(w(X))$. If $\beta(w(X))=0$ this coefficient is zero. If $\beta(w(X))>0$ the coefficient is positive, and if $\beta(w(X))<0$ it is negative. Therefore, basis vectors of $\k$ which are $(+1)$-eigenvectors of $\theta\tau$ contribute the set $\{ -\beta(w(X)) H_\beta : \beta \in \Delta^+(\g,\a^{-\tau}), (\g^\beta)_+ \neq 0 \}$ to $\Gamma_{w(X)}$.

The fact that $\langle pr_\h [Z,Ad(w(a))v_i],Ad(w(a))v_j \rangle = 0$ holds for all $Z\in \t $ whenever $i\neq j$ follows from general properties of the Killing form.
 
The conclusion is that the cone $\Gamma_{w(X)}=\Phi(V_{w(a).K}^*)$ is generated by the weights
\begin{eqnarray*}
\{ -\beta(w(X)) H_\beta : \beta \in \Delta^+(\g,\a^{-\tau}), (\g^\beta)_+ \neq 0 \} \\
\cup \ \{ -H_{\beta}: \beta \in \Delta^+(\g,\a^{-\tau}), (\g^\beta)_- \neq 0 \} ,
\end{eqnarray*}
as asserted.

\end{proof}

\begin{cor}\label{conescor}
The image of the moment map $\Phi(M_a)$ is contained in the set $w'(X)+\Gamma_+$, where $w'\in {\mathcal W}_{K\cap H}$ is such that $\beta(w'(X))\geq 0$ for all $\beta \in \Delta^+(\g,\a^{-\tau})$ and $\Gamma_+$ is the proper cone $\Gamma_+=\rm{cone}(-H_\beta : \beta \in \Delta^+(\g,\a^{-\tau}))$.
\end{cor}
\begin{proof}
From Theorem \ref{hnp} and Remark \ref{hnpremark} we know that there is a neighborhood $U_{w'(X)}\subseteq \a^{-\tau}$ of $w'(X)$ such that $\Phi(M_a)\cap U_{w'(X)}= (w'(X)+\Gamma_{w'(X)}) \cap U_{w'(X)}$. Lemma \ref{cones} implies that $\Gamma_{w'(X)} \subseteq \Gamma_+$. Suppose there exists some $Z\in \Phi(M_a)$ such that $Z\not\in w'(X)+\Gamma_+$. Since $\Phi(M_a)$ is convex the line segment $\overline{w'(X)Z}$ lies entirely in $\Phi(M_a)$. Fix some $Y\in \overline{w'(X)Z}\cap U_{w'(X)}$ with $Y\neq w'(X)$. Then $Y\in \Phi(M_a)\cap U_{w'(X)} \subseteq w'(X)+\Gamma_{w'(X)} \subseteq w'(X)+\Gamma_+$. But this implies $Z\in w'(X)+\Gamma_+$ since $\Gamma_+$ is a cone and $Y\neq w'(X)$, a contradiction. Therefore, $\Phi(M_a)\subseteq w'(X)+\Gamma_+$. The cone $\Gamma_+$ is proper since it is spanned by  vectors $-H_\beta$ associated to positive roots $\beta$.
\end{proof}

The special property of $\Phi(M_a)$ stated in the corollary allows us to describe $\Phi(M_a)$ entirely in terms of the local cones $\Gamma_{w(X)}$ associated to the fixed points, as the following proposition shows.

\begin{prop}\label{semibounded}
Let $C$ be a closed, convex, locally polyhedral set (in some finite dimensional vector space $V$). Denote by $\Gamma_c$ the local cone at $c\in C$ ( i.e. there is a neighborhood $U_c\subset V$ of $c$ such that $C\cap U_c = (c+\Gamma_c)\cap U_c $). Suppose $C\subset x+\Gamma$ for some $x\in V$ and some proper cone $\Gamma \subset V$. Then
\[ C=\bigcap_{\Gamma_c \ \rm{proper}} (c+\Gamma_c) , \]
i.e. $C$ is completely determined by the local cones that are proper.
\end{prop}

\begin{proof}
For any $c\in C$ we write $d_c$ for the dimension of the maximal subspace contained in $\Gamma_c$. (In particular, $d_c=0$ means that $\Gamma_c$ is proper.) First we will show that if $d_c>0$, then $c\in c'+\Gamma_{c'}$ for some $c'$ with $d_{c'}<d_c$.

If $d_c>0$, then $\Gamma_c$ contains a line, say $L$. Since $C$ lies in a proper cone, $(c+L)\cap C$ is semi-bounded. We pick an endpoint $c'$ of $(c+L)\cap C$. Since $C$ is closed $c'\in C$, and clearly $c\in c'+\Gamma_{c'}$. Convexity of $C$ implies that if a line $L'$ is contained in $\Gamma_{c'}$ then $L'\subset \Gamma_{\tilde{c}}$ for each inner point $\tilde{c}$ of $(c+L)\cap C$. In particular, $d_{c'}\leq d_c$. On the other hand, $\Gamma_{c'}$ does not contain the line $L\subset \Gamma_c$. Therefore, $d_{c'}<d_c$.

Now, the assumptions on $C$ imply
\[ C=\bigcap_{c\in C} (c+\Gamma_c) \]
If we set $n=\rm{dim}(V)$ the above arguments lead to
\[ C=\bigcap_{d_c\leq n} (c+\Gamma_c)=\bigcap_{d_c\leq n-1} (c+\Gamma_c) = \dots = \bigcap_{d_c=0} (c+\Gamma_c) \]
\end{proof}

We are now ready to give the desired description of $\Phi(M_a)$ which is the content of van den Ban's theorem.

\begin{theorem}
The set $\Phi(M_a)=(pr_{\a^{-\tau}} \circ \log \circ \mu)(Ha)$ is the sum of a compact convex set and a closed (proper) cone $\Gamma$. More precisely, for $a=\exp X$,
\[ \Phi(M_a)={\rm conv}({\mathcal W}_{K\cap H}.X)+\Gamma , \]
with
\[ \Gamma= {\rm cone} \{ -H_\beta : \beta \in \Delta^+(\g,\a^{-\tau}), (\g^\beta)_- \neq 0 \} \]
\end{theorem}

\begin{proof}
The image $\Phi(M_a)$ is closed, convex and locally polyhedral. Moreover, by Corollary \ref{conescor}, it is contained in $w'(X)+\Gamma_+$ for some proper cone $\Gamma_+$. Proposition \ref{semibounded} implies that $\Phi(M_a)$ is determined by the local cones that are proper. According to Remark \ref{hnpremark}, a local cone $\Gamma_{\Phi(m)}$ can be proper only if $\t_m=\t$, i.e. if $m$ is a $T$-fixed point. The $T$-fixed points have been characterized in Proposition \ref{fixed}, so Proposition \ref{semibounded} yields
\[ \Phi(M_a)=\bigcap_{w\in {\mathcal W}_{K\cap H}} (w(X)+\Gamma_{w(X)}), \]
with $\Gamma_{w(X)}$ as in Lemma \ref{cones}. 

The sum ${\rm conv}({\mathcal W}_{K\cap H}.X)+\Gamma$ is closed, convex and locally polyhedral as well. As a sum of a compact set and the proper cone $\Gamma$ it is contained in $x+\Gamma$ for some $x\in \a^{-\tau}$, hence Proposition \ref{semibounded} is applicable. First we want to see at which points in ${\rm conv}({\mathcal W}_{K\cap H}.X)+\Gamma$ the local cone is proper. Let $c\in {\rm conv}({\mathcal W}_{K\cap H}.X)$ and $\gamma \in \Gamma$. Clearly, the local cone at $c+\gamma$ is improper unless $\gamma=0$. But then $c+\gamma=c$ is contained in a convex set with extremal points $\{ w(X) : w\in {\mathcal W}_{K\cap H} \} $. The local cone can be proper only if $c+\gamma$ is one of those extremal points. Proposition \ref{semibounded} now gives
\[ {\rm conv}({\mathcal W}_{K\cap H}.X)+\Gamma = \bigcap_{w\in {\mathcal W}_{K\cap H}} (w(X)+\Gamma'_{w(X)}). \]
Here, $\Gamma'_{w(X)}$ denotes the local cone of ${\rm conv}({\mathcal W}_{K\cap H}.X)+\Gamma$ at $w(X)$. To finish the proof it is sufficient to show that $\Gamma'_{w(X)}=\Gamma_{w(X)}$.

Clearly, $\Gamma'_{w(X)}=\Gamma''_{w(X)}+\Gamma$, where $\Gamma={\rm cone} \{ -H_\beta : \beta \in \Delta^+(\g,\a^{-\tau}), (\g^\beta)_- \neq 0 \}$ as before and $\Gamma''_{w(X)}= {\rm cone} \{ w'(X)-w(X): w'\in {\mathcal W}_{K\cap H} \} $. From Lemma \ref{cones} we know that $\Gamma_{w(X)}$ contains the cone $\Gamma$. Moreover, the set $\Phi(M_a)$ is convex and contains all points $w(X)$, and therefore contains ${\rm conv}({\mathcal W}_{K\cap H}.X)$. This implies that its local cone at $w(X)$, i.e. $\Gamma_{w(X)}$, contains $\Gamma''_{w(X)}$ as well. Therefore, $\Gamma_{w(X)}\supseteq \Gamma''_{w(X)}+\Gamma=\Gamma'_{w(X)}$.

Each root $\beta\in \Delta(\g,\a^{-\tau})$ defines the isomorphism 
\[ s_\beta:\a^{-\tau} \rightarrow \a^{-\tau}, Z\mapsto Z-2\frac{\beta{Z}}{\langle \beta,\beta \rangle}H_\beta . \]
In view of Remark \ref{thetatau} the Weyl group ${\mathcal W'}={\mathcal W}_{K\cap H}$ consists exactly of those $s_\beta$ for which $(\g^\beta)_+ \neq 0$. In particular, $s_\beta(w(X))\in {\mathcal W}_{K\cap H}$ for all $\beta\in \Delta^+(\g,\a^{-\tau})$ for which $(\g^\beta)_+ \neq 0$. The identity $s_\beta(w(X))-w(X)=-2\frac{\beta(w(X))}{\langle \beta,\beta \rangle}H_\beta$ implies ${\rm cone} \{ -\beta(w(X))H_\beta: \beta \in \Delta^+(\g,\a^{-\tau}), (\g^\beta)_+ \neq 0 \} \subseteq \Gamma''_X$. With Lemma \ref{cones} we obtain $\Gamma_{w(X)}\subseteq \Gamma''_{w(X)}+\Gamma=\Gamma'_{w(X)}$.
\end{proof}

\bibliographystyle{amsplain}

\Addresses
\end{document}